\documentclass[10pt]{article}

\usepackage{a4wide}
\usepackage{amssymb}
\usepackage{amsfonts}
\usepackage{amsmath}
\input xy
\xyoption{arrow} \xyoption{matrix}

\date{}

\newtheorem{proposition}{Proposition}[section]
\newtheorem{theorem}[proposition]{Theorem}
\newtheorem{lemma}[proposition]{Lemma}

\newtheorem{corollary}[proposition]{Corollary}

\def\der{\partial }

\def\nFM0{{\nu }_{F,M_0}}
\def\nFN0{{\nu }_{F,N_0}}
\def\nGN0{{\nu }_{G,N_0}}

\def\N0{ {\bf N}_0 }

\def\g{\gamma}
\def\v{\varphi}
\def\ra{\rightarrow}

\def\lra{\leftrightarrow}
\def\Xpm{X^{\pm }}

\def\s{\sigma}
\def\Z{\mathbb{Z}}

\def\l1{{\lambda}_1}

\def\a{\alpha}
\def\a0{ {\alpha }_0}
\def\a1{ {\alpha }_1}

\def\l{\lambda}

\def\nFGM0{{\nu }_{F,G,M_0}}

\def\nFN0{{\nu}_{F,N_0}}

\def\sm{{\sigma}^m}

\def\sm1{{\sigma}^{-1}}

\def\smtp1{{\sigma}^{-t+1}}

\def\S1{S^{-1}}

\def\Xpm1{X^{\pm 1}_1}

\def\sPM1{{\sigma }^{\pm 1}}
\def\sMP1{{\sigma }^{\mp 1 }}

\def\d{\delta}

\def\di{{\rm d.ind}}

\def\L{\Lambda}

\def\Ytm1{Y^{t-1}}
\def\Yim1{Y^{i-1}}

\def\CF{{\cal F}}

\def\supp{{\rm supp }}

\def\Der{{\rm Der }}
\def\ad{{\rm ad }}
\def\dim{{\rm dim }}
\def\char{{\rm char }}

\def\ker{ {\rm ker } }

\def\gcd{ {\rm gcd } }
\def\D{ \Delta }
\def\Ev{ {\rm Ev} }

\def\SL2Z{ {\rm SL}_2({\bf Z}) }

\def\th{ \theta }

\def\Gp1{ G^{1 , 1 } }
\def\P11{ P^{-1 , 1 } }
\def\Pp1{ P^{1 , 1 } }

\def\NP{{\rm NP}}

\def\th{\theta}

\def\nCLsr{{}^\nu\kern-2pt {\cal L}^{\sigma , \rho  }}
\def\nP{{}^\nu \kern-2pt P}
\def\nL{{}^\nu\kern-2pt L}
\def\nLL{{}^\nu\kern-2pt \Lambda}
\def\nPsr{{}^\nu\kern-2pt P^{\sigma , \rho  }}
\def\nLsr{{}^\nu\kern-2pt L^{\sigma , \rho  }}
\def\nuCL{{}^\nu\kern-2pt  {\cal L}}
\def\nCLsr{{}^\nu\kern-2pt {\cal L}^{\sigma , \rho  }}
\def\nCL1m{{}^\nu\kern-2pt {\cal L}^{-1 , 1  }}
\def\x1nu{x^\frac{1}{\nu}}
\def\xm1nu{x^{-\frac{1}{\nu}}}

\def\ra{\rightarrow }

\def\CB{{\cal B}}

\def\coker{{\rm coker}}

\def\nAM0{{\nu }_{{\cal A},M_0}}
\def\nAN0{{\nu }_{{\cal A},N_0}}

\def\Der{ {\rm Der }}

\def\ad{ {\rm ad }}

\def\SL{{\rm SL}}

\def\di!{\frac{\der^i}{i!}}
\def\dik!{\frac{\der^k_i}{k!}}

\def\N{\mathbb{N}}

\def\0{\overline{0}}
\def\1{\overline{1}}

\def\Ln1{\L_{n,\overline{1}}}

\def\a1{a_{\overline{1}}}

\def\S{\Sigma}

\def\vn1{\overrightarrow{n-1}}

\def\im{{\rm im}}

\def\mJ{\mathbb{J}}
\def\mI{\mathbb{I}}

\def\ind{{\rm ind}}

\def\Frac{{\rm Frac}}

\def\K1{{\rm K}_1}

\def\hmI1{\widehat{\mI_1}}
\def\tmI1{\widetilde{\mI_1}}
\def\tmJ1{\widetilde{\mJ_1}}
\def\hB1{\widehat{B_1}}
\def\hCB1{\widehat{\CB_1}}

\begin{document}

\author{V. V. \  Bavula 
}

\title{On the eigenvector algebra  of the product of elements with commutator one in the first Weyl algebra }

\maketitle

\begin{abstract}
Let $A_1=K \langle X, Y \, | \, [Y,X]=1\rangle$ be the (first)
Weyl algebra over a field $K$ of characteristic zero. It is known
that the set of eigenvalues of the inner derivation $\ad (YX)$  of
$A_1$  is $\Z$. Let $ A_1\ra A_1$, $X\mapsto x$, $Y\mapsto y$, be
a $K$-algebra homomorphism, i.e. $[y,x]=1$. It is proved that  the
set of eigenvalues of the inner derivation $\ad (yx)$ of the Weyl
algebra $A_1$ is $\Z$  and the eigenvector algebra of $\ad (yx)$
is $K\langle x,y\rangle $ (this would be an easy corollary of the
Problem/Conjecture of Dixmier of 1968 [still open]: {\em is an
algebra endomorphism of $A_1$ an automorphism?}).




\end{abstract}


\section{Introduction}

The following notation is fixed throughout the paper:
 $A_1=K \langle X, Y \, | \, [Y,X]=1\rangle$ is  the (first) {\em
Weyl algebra} over a field $K$ of characteristic zero where
$[a,b]:= ab -ba$ is the {\em commutator} of  elements $a$ and $b$;
$H:=YX$;  $K^*:=K\backslash \{ 0\}$;
 $C(a):= \{ b\in A_1\, | \, ab=ba\}$ is
the {\em centralizer} of $a$ in $A_1$; $\ad (a) := [a, \cdot ]$ is
the {\em inner derivation}  associated with an element $a$; the
subalgebra $D(a)$ of $A_1$ generated by all the eigenvectors of
$\ad (a)$ in $A_1$ is called the {\em eigenvector algebra} of $\ad
(a)$ or $a$;  a pair of elements $x$ and $y$ of the Weyl algebra
$A_1$ such that $[y,x]=1$ or, equivalently, a $K$-algebra
endomorphism $A_1\ra A_1$, $X\mapsto x$, $Y\mapsto y$; $A_1':=
K\langle x,y \rangle$ is the Weyl algebra and $\Frac (A_1')$ is
its total skew field of fractions ($A_1'\subseteq A_1$ and $\Frac
(A_1')\subseteq \Frac (A_1)$).

$\noindent $

{\bf The First Problem/Conjecture of Dixmier, \cite{Dix}}: {\em is
an algebra endomorphism of the Weyl algebra $A_1$ an
automorphism?}

$\noindent $

For a connection of the Conjecture of Dixmier with the Jacobian
Conjecture the reader is referred to the papers of Tsuchimoto
\cite{Tsuchi05}, Belov-Kanel and Kontsevich \cite{Bel-Kon05JCDP},
and the author \cite{JC-DP}.

$\noindent $

The aim of the paper is to prove the following four theorems.

\begin{theorem}\label{18Dec10}
 If $[y,x]=1$ for some elements $x$ and $y$ of the Weyl
algebra $A_1$ then $C(yx)=K[yx]$.
\end{theorem}

\begin{theorem}\label{22Dec10}
 If $[y,x]=1$ for some elements $x$ and $y$ of the Weyl
algebra $A_1$ then eigenvector algebra  $D(yx)$ is equal to the
subalgebra $K\langle x,y\rangle$ of $A_1$.
 In particular, the set of eigenvalues of the inner derivation
 $\ad (yx)$ of the Weyl algebra $A_1$  is $\Z$, and, for each $i\in \Z$, the vector space of
 eigenvectors for $\ad (yx)$ with eigenvalue $i$ is $K[yx]v_i'$
 where $v_i':=x^i$ if $i\geq 0$, and
$v_i:=y^{-i}$ $if <0$.
\end{theorem}

\begin{theorem}\label{X1Jan11}
$\Frac (A_1') \cap A_1= A_1'$ where the intersection  is taken in
$\Frac (A_1)$.
\end{theorem}

Dixmier \cite{Dix} proved that each maximal commutative subalgebra
$C$ of the Weyl algebra $A_1$ coincides with the centralizer
$C(a)$ of every non-scalar element $a$ of $C$. Let $V$ be a vector
space. A linear map $\v : V\ra V$ is called a {\em locally
nilpotent map} if $V=\bigcup_{n\geq 1}\ker (\v^n)$, i.e. for each
element $v\in V$ there exists a natural number $n$ such $\v^nv=0$.
 The linear map $\v $ is called a {\em semi-simple} linear map if
 $V= \bigoplus_{\l \in \Ev (\v )}\ker (\v -\l )$ where $\Ev (\v )$
  is the set of eigenvalues of the map $\v $ in the field $K$. The
  next theorem classifies (up to isomorphism) the maximal
  commutative subalgebras $C$ of the Weyl algebra $A_1$ that admit
  either a locally nilpotent derivation or a semi-simple
  derivation $\d$ (or both). In each case, the derivation $\d$
  and its eigenvalues are found.

\begin{theorem}\label{23Dec10}
Let $C$ be a maximal commutative subalgebra of the Weyl algebra
$A_1$ and $0\neq \d \in \Der_K(C)$. Then
\begin{enumerate}
\item The derivation $\d$ is a locally nilpotent derivation of the
algebra $C$ iff there exists an element $c\in C\backslash \ker (\d
)$ such that $\d^n (c)=0$ for some $n\geq 2$ iff $C = K[t]$ for
some element $t\in C$ and $\d = \frac{d}{dt}$. \item The
derivation $\d$ is a semi-simple derivation of the algebra $C$ iff
$\d (c) = \l c$ for some $0\neq c\in C$ and $\l \in K^*$ iff
$C=\bigoplus_{\l \in \Ev (\d )} \ker (\d -\l )$, $\dim_K(\ker (\d
-\l ))= 1$ for all $\l \in \Ev (\d )$, and the additive monoid
$\Ev (\d )$ of eigenvalues of the derivation $\d$ is a submonoid
of the additive monoid $\N \rho$ for some $\rho \in K^*$ such that
$\Z \Ev (\d ) = \Z \rho $.  Moreover, one of the following three
cases occurs (the cases (ii) and (iii) are not mutually
exclusive):

(i) if $C= K[H]$ where $H:=YX$  then $\d = \rho H \frac{d}{dH}$
and $\Ev (\d ) = \rho \N$;

(ii) if $0\varsubsetneqq \supp (C) \subseteq \N$ then $\Ev (\d ) =
\rho' \supp (C)$ for some $\rho'\in K^*$. Moreover, $\Ev ( \d ) =
\frac{\l }{v (c_\l )}\supp (C)$ for all $0\neq \l \in \Ev ( \d )$
and $ 0\neq c_\l \in \ker (\d - \l )$; and $\d $ is the unique
extension of the derivation $\l c_\l \frac{d}{dc_\l}$ of the
polynomial algebra $K[c_\l]$ to $C$;

(iii) If $0\varsubsetneqq \supp_-(C)\subseteq -\N$ then $\Ev (\d )
= \rho' \supp_- (C)$ for some $\rho'\in K^*$. Moreover, $\Ev ( \d
) = \frac{\l }{v_- (c_\l )}\supp_- (C)$ for all $0\neq \l \in \Ev
( \d )$ and $ 0\neq c_\l \in \ker (\d - \l )$; and $\d $ is the
unique extension of the derivation $\l c_\l \frac{d}{dc_\l}$ of
the polynomial algebra $K[c_\l]$ to $C$;

iff there exists a nonzero additive submonoid $E$ of $(\N , +)$
such that the algebra $C$ is isomorphic to the monoid subalgebra
$\bigoplus_{i\in E}Kt^i$ of the polynomial algebra $K[t]$  in a
variable $t$, $\d = \rho t\frac{d}{dt}$ for some $\rho \in K^*$
and $\Ev (\d ) = \rho E$. \item $\ker (\d ) =K$.
\end{enumerate}
\end{theorem}

{\it Remarks}. 1.  In general, it is not true that in the cases
2(ii),(iii) the algebra $C$ is isomorphic to a polynomial algebra
or is a non-singular algebra (see Example 4, Section
\ref{DROPFUN}).

2. Theorem \ref{18Dec10} is a particular case  Theorem 1.2,
\cite{joseph-AJM}. Here a shorter and different proof is given.

As a corollary  two (short) proofs are given to the following
theorem of J. A. Guccione, J. J. Guccione and C. Valqui.

\begin{theorem}\label{GGV-Cen}
\cite{GGV} If $[y,x]=1$ for some elements $x$ and $y$ of the Weyl
algebra $A_1$ then $C(x)=K[x]$.
\end{theorem}

If the Problem/Conjecture  of  Dixmier  is true then Theorems
\ref{18Dec10}, \ref{22Dec10}, \ref{X1Jan11} and \ref{GGV-Cen}
would be its easy corollaries.

$\noindent $

{\bf Proof I of Theorem \ref{GGV-Cen}}. The inner derivations $\ad
(x)$ and $\ad (y)$ of the Weyl algebra $A_1$  commute: $0=\ad (1)
= \ad ([y,x])=[\ad (y), \ad (x)]$. So, $\ad (y)$ is a non-zero
derivation of the algebra $C(x) = \ker ( \ad (x)) $ such that $\ad
(y)^2(x)=[y,[y,x]]=[y,1]=0$ and $x\in C(x) \backslash \ker ( \ad
(y))$ since $\ad (y) (x)=1$. By Theorem \ref{23Dec10}.(1), $C(x) =
K[t]$ for some element $t\in C(x)$ such that $ [y,t]=1$. Then $[y,
t-x]=0$, and so $t-x \in K$, by Theorem \ref{23Dec10}.(3) or by
Theorem \ref{2facts}.(2) (since $t-x\in C(x) \cap C(y) =K$), and
so $C(x) = K[x]$. $\Box$

$\noindent $

{\bf Proof II of Theorem \ref{GGV-Cen}}. The proof follows from
the following two facts:

\begin{theorem}\label{2facts}
\begin{enumerate}
\item (Lemma 2.1, \cite{invfor-JPAA-2007}) Let $R$ be a ring $R$
and $\d$ be a locally nilpotent derivation such that $\d (x)=1$
for some element $x\in R$. Then the ring $R$ is the skew
polynomial ring $\ker (\d ) [ x; \ad (x)]$. \item \cite{Dix} Let
$a$ and $b$ be non-scalar elements of the Weyl algebra $A_1$ such
that $ab \neq ba$. Then $C(a) \bigcap C(b) = K$.
\end{enumerate}
\end{theorem}
Since $[y,x]=1\neq 0$, $C(x) \bigcap C(y) = K$ (by the second
fact) and then, by the first fact where $R=C(x)$ and $\d = \ad
(y)$ ($\d$ is a locally nilpotent derivation of $R$  by Theorem
\ref{23Dec10} since $\d^2(x)=0$ and $\d (x) =1\neq 0$):
$$ C(x) = C(x)\cap C(y) [ x; \ad (x)=0] = K[x].\;\;\; \Box$$
All the results of Sections \ref{DROPFUN} and \ref{PR22DEC} can be
seen as intermediates steps in the proofs of Theorems
\ref{18Dec10}, \ref{22Dec10}, \ref{X1Jan11} and \ref{23Dec10}. One
of the key ideas in the proofs is to show that certain maps are
locally nilpotent. The concept of the {\em drop} of a linear map
and the linear maps of {\em constant drop} are an important tool
in proving that certain linear maps are locally nilpotent (see
Theorem \ref{20Dec10}). It turns out that the drop of a linear map
coincides  up to a negative rational multiple with the index of
the map (Theorem \ref{20Dec10}.(2)). Theorems \ref{18Dec10} and
\ref{GGV-Cen} beg to ask the following question (see also
Questions 2 and 3 at the end of the paper).

$\noindent $

{\it Question 1. Is it true that if $ab = ba$ for some elements
$a\in A_1'$ and $b\in A_1$ then $b\in A_1'$, i.e. $C_{A_1'}(a) =
C_{A_1}(a)$? }


\section{Linear maps of constant drop, proof of Theorem \ref{18Dec10}}\label{DROPFUN}

In this section proofs of Theorems \ref{18Dec10} and \ref{23Dec10}
are given. Many results of this section are used in the proof of
Theorem \ref{22Dec10} which is given in Section \ref{PR22DEC}.

Let a $K$-algebra $A$ be a domain (not necessarily commutative). A
function $v: A\backslash \{ 0\} \ra \Z$ is called a {\em degree
function} on $A$, if for all elements $a,b\in A\backslash \{ 0\}$,
\begin{enumerate}
\item  $v(ab) = v (a) + v(b)$, \item $v(a+b) \leq \max \{ v(a) ,
v(b)\}$ where $v(0):=-\infty$, and \item  $v(\l ) =0$ for all $\l
\in K^*:=K\backslash \{ 0\}$.
\end{enumerate}
The degree function $v$ on $A$ determines the ascending algebra
filtration $A= \bigcup_{i\in \Z} A_{\leq i}$ which is called the
$v$-{\em filtration} ($A_{\leq i}A_{\leq j} \subseteq A_{\leq
i+j}$ for all $i,j\in \Z$) where $A_{\leq i} :=\{ a\in A\, | \,
v(a) \leq i\}$,  $\bigcap_{i\in \Z} A_{\leq i}=0$ and the
associated graded algebra ${\rm gr}_v (A) := \bigoplus_{i\in \Z}
{\rm gr}_{v,i}(A)$ is a domain where ${\rm gr}_{v, i}(A) :=A_{\leq
i } / A_{\leq i-1}$. If $v(a)\neq v(b)$ then $v(a+b) = \max \{
v(a) , v(b)\}$.

$\noindent $

{\it Example 1}. The usual degree $\deg_x$ of a polynomial is a
degree function on the polynomial algebra $K[x]$ or $K[x=x_1, x_2,
\ldots , x_n]$. The {\em total degree} $\deg$ on $K[x_1, x_2,
\ldots , x_n]$ is another example of a degree function. More
generally, for each nonzero vector $d= (d_1, \ldots , d_n) \in
\Z^n$, the function $v_d: K[x_1, x_2, \ldots , x_n]\ra \Z $,
$$ v_d(\sum \l_{\alpha_1, \ldots ,\alpha_n}x_1^{\alpha_1}\cdots
x_n^{\alpha_n}):=\max \{ d_1\alpha_1+\cdots + d_n\alpha_n\, | \,
\l_{\alpha_1, \ldots ,\alpha_n} \neq 0\}$$ is a degree function
where $\l_{\alpha_1, \ldots ,\alpha_n}\in K$.

 $\noindent $

{\it Example 2}. Let $\nu$ be a {\em discrete valuation} on the
 $K$-algebra $A$ which is a domain, i.e. $\nu : A\backslash \{
 0\}\ra \Z$ is a function such that for all elements $a,b\in A\backslash \{ 0\}$,
\begin{enumerate}
\item  $\nu (ab) = \nu (a) + \nu(b)$, \item $\nu (a+b) \geq \min
\{ \nu (a) , \nu (b)\}$ where $\nu (0):=\infty$, and \item $\nu
(\l ) =0$ for all $\l \in K^*$.
\end{enumerate}
Then $v:= -\nu$ is a degree function on $A$, and vice versa.

$\noindent $

{\it Example 3}. Let a $K$-algebra $A=\bigoplus_{i\in \Z} A_i$ be
a $\Z$-graded domain with $ K\subseteq A_0$. Then the $\Z$-{\em
graded degrees} $v, v_-: A\ra \Z$,  
\begin{equation}\label{vamin}
v(a) := \min \{ d\in \Z \, | \, a\in \bigoplus_{i\leq d} A_i\},
\;\;\; v_-(a) := -\max \{ d\in \Z \, | \, a\in \bigoplus_{i\geq d}
A_i\}
\end{equation}
are  degree functions on $A$. In particular, we have the
$\Z$-graded degree $v$ on the Weyl algebra $A_1$ (see below). Let
$C$ be a subalgebra of the algebra $A$. Then $\supp (C):= v
(C\backslash \{ 0\} )$ and $\supp_-(C):= - v_- (C\backslash \{ 0\}
)$ are additive submonoids of $(\Z , +)$.

$\noindent $

Let $C$ be a $K$-subalgebra of the algebra $A$. For a $K$-linear
map $\d : C\ra C$, define the map $\D_{\d , v}$,  which is called
the {\em drop} of $\d$ with respect to the degree function $v$ on
the algebra $A$, by the rule 
\begin{equation}\label{ddef}
\D = \D_{\d , v} : C\ra \Z \bigcup \{ -\infty \}, \;\; a\mapsto v(
\d (a)) - v(a),\;\; \D (0):= -\infty ,
\end{equation}
i.e. $v(\d (a))= v(a)+ \D (a)$ for all elements $a\in C$ where
$-\infty +z = -\infty$ for all elements $z\in \Z$.  For all
elements $a\in C$ and $\l \in K^*$, $\D (\l a) = \D (a)$. For all
elements $a,b\in C$ such that either $v(a)\neq v(b)$ or $v(a) =
v(b) = v(a+b)$, 
\begin{equation}\label{daba}
\D (a+b) \leq \max \{ \D (a), \D (b)\}.
\end{equation}
For linear maps $\d , \d' : C\ra C$ and for all elements $a\in C$,
\begin{equation}\label{daba1}
\D_{\d \d', v}(a) =  \D_{\d, v}(\d' (a)) + \D_{\d', v}(a).
\end{equation}

\begin{lemma}\label{a17Dec10}
Let $A$, $C$, $v$, $\d$ and $\D =\D_{\d , v}$ be as above. Suppose
that $\d$ is a derivation of the algebra $C$. Then
\begin{enumerate}
\item $\D (ab) \leq \max \{ \D (a) , \D (b) \} $ for all elements
$a, b\in A$.  \item If, in addition, $C$ is a commutative algebra
and $\char (K)=0$ then $\D (a^n) = \D (a)$ for all elements $a\in
C$ and $n\geq 1$.
\end{enumerate}
\end{lemma}

{\it Proof}. 1. \begin{eqnarray*}
 \d  (ab) &=& v ( \d (ab))-v(ab) = v( \d (a) b +a\d (b)) - v(a) - v(b)
 \leq \max \{ v( \d (a) b), v(a\d (b))\}  - v(a) - v(b)\\
 &=& \max \{  v( \d (a))+v( b) - v(a) - v(b), v(a)+v(\d (b)) - v(a) -
 v(b)\}\\
 &= & \max \{ \D (a), \D (b)\}.
\end{eqnarray*}
2.  $\d  (a^n) = v ( \d (a^n))-v(a^n) = v( na^{n-1}\d (a)) - nv(a)
= (n-1) v(a) +v(\d (a)) - n v(a) = \D (a)$.  $\Box $

$\noindent $

{\it Definition}. The linear map $\d : C\ra C$ is called a {\em
linear map of constant drop} if $\D (a) = \D (b)$ for all elements
$a,b\in C\backslash \ker (\d )$, and the common value of $\D (a)$
where $a\in C \backslash \ker (\d )$ is called the {\em drop} of
the linear map $\d$ denoted by $\D = \D_{\d , v}$.

$\noindent $

Let $V$ be a vector space over the field $K$, a linear map $\v :
V\ra V$ is called a {\em Fredholm map/operator} if it has finite
dimensional kernel and cokernel, and $$\ind (\v ) := \dim_K( \ker
(\v )) - \dim_K(\coker (\v ))$$ is called the {\em index} of the
map $\v$. Let $\CF (V)$ be the set of all Fredholm linear maps in
$V$, it  is, in fact, a monoid since 
\begin{equation}\label{indpps}
\ind (\v \psi ) = \ind (\v ) +\ind (\psi )\;\; {\rm for \; all}
\;\; \v , \psi \in \CF (V).
\end{equation}
The next theorem provides examples of Fredholm  linear maps of
constant drop and of locally nilpotent maps.

\begin{theorem}\label{20Dec10}
Let $A$, $C$, $v$, $\d$ and $\D_{\d , v}$ be as above and $C':=
C\backslash K$. Suppose that the following conditions hold.
\begin{enumerate}
\item $\D (ab) \leq \max \{ \D (a), \D (b)\} $ for all elements
$a, b\in C$. \item $\D (a^n) = \D (a)$ for all elements $a\in C$
and $n\geq 1$. \item $v(c)\geq 0$ for all elements $c\in C$, and
$v(C') \neq 0$. \item For all $i\geq 0$, $\dim_K(C_{\leq i} /
C_{\leq i-1})\leq 1$ where $C_{\leq i}:= C\cap A_{\leq i} = \{
c\in C\, | \, v(c) \leq i\}$. \item $\ker (\d ) =K$.
\end{enumerate}
Then
\begin{enumerate}
\item  $\D (a) = \D (b)$ for all elements $a, b\in C'$, i.e. the
map $\d$ has constant drop $\D = \im (C')$. \item The map $\d \in
\CF (C)$ is a Fredholm map with index $\ind (\d ) = -\frac{\D
}{g}$ and $\dim_K( \coker (\d )) = \frac{\D}{g}+1$ where $g$ is
the unique positive integer such that $\Z v(C') = \Z g$. Moreover,
there exists an explicit natural number $\g$ (see (\ref{gmud}))
and a vector subspace  $V$ of $C_{\leq \g }$ such that $C=
V\bigoplus \im (\d )$, i.e. $V\simeq \coker (\d )$. \item If $\D
(a)<0$ for some element $a\in C$, i.e. $\D <0$, then $\D = -g$ and
the map $\d$ is a locally nilpotent map such that there exists a
$K$-basis $\{ e_i\}_{i\in \N}$ of the algebra $C$ such that, for
all $i\geq 0$, $\d (e_i) = e_{i-1}$ and $C_{\leq ig } =
\bigoplus_{j=0}^iKe_j$ where $e_{-1}:=0$. \item {\rm (Additivity
of the drop)} Let $ \d': C\ra C$ be another linear map that
satisfies conditions 1--5 then $\D_{\d \d', v} = \D_{\d , v} +
\D_{\d', v}$, i.e.  for all elements $a\in C\backslash \ker (\d
\d')$, $v( \d\d' (a))= v(a) +\D_{\d , v} + \D_{\d', v}$.
\end{enumerate}
\end{theorem}

{\it Proof}. 1.  By condition 3, $C\cap A_{\leq -1}=0$. Then, by
condition 4, $C\cap A_{\leq 0} =K$, and so the map $v: C\backslash
\{ 0\} \ra \N$, $ c \mapsto v(c)$, is a non-zero homomorphism of
monoids (by condition 3). Let $H$ be its image. The $\Z$-submodule
$\Z H$ of $\Z$ is equal to $\Z g$ where $g=\gcd \{ i \, | \, i\in
H\}$. It is a well-known fact (and easy to show)  that
\begin{equation}\label{gNH}
| g\N \backslash H | <\infty.
\end{equation}
In more detail, let $g=s-t$ for some elements $s,t\in H$. Let $m$
be the least positive element of $H$. For each $k=0, 1, \ldots ,
mg^{-1}-1$, the element $h_k:= t (mg^{-1} -k) +ks\in H$. Since
$h_k = tg^{-1} m + k(s-t) = tg^{-1} m + kg\equiv kg \mod m$ for
all $k=0, 1, \ldots , mg^{-1}-1$, it follows that $|g\N \backslash
\bigcup_{k=0}^{mg^{-1}-1} (h_k+\N m )|<\infty $,  and (\ref{gNH})
follows since 
\begin{equation}\label{gNH1}
H':=\bigcup_{k=0}^{mg^{-1}-1} (h_k+\N m )\subseteq H.
\end{equation}
Claim: {\em there  exists an integer $l$ such that $\D (c) \leq l$
for all $c\in C$}.

$\noindent $

To prove this fact choose elements $a,b\in C$ such that $v(a) = m$
and $v(b) +m\Z$ is a generator for the finite group $\Z H / \Z m=
\Z g / \Z m$. Let $ g_1= v(b)$. Then $|H\backslash
\bigcup_{i=0}^{mg^{-1}-1} (ig_1+\N m )|<\infty $, and so $\dim_K(
C/ \bigoplus_{j\in \N } \bigoplus_{i=0}^{mg^{-1}-1}
Ka^jb^i)<\infty$, and so  the algebra $C$  has a $K$-basis of the
type $\{ e_1, \ldots , e_t, a^jb^i\, | \, j\in \N , i=0, 1, \ldots
, mg^{-1}-1\}$ where the degree function $v$ takes {\em distinct}
values on the elements of the basis. By (\ref{daba}) and
conditions 1 and 2, for all elements $c\in C'$,
$$ \D (c) \leq l:=\max \{ \D(e_1), \ldots , \D (e_t), \D (a), \D
(b)\}.$$ Fix an element $c\in C'$ such that $\D (c)$ is the
largest possible. Then $p:=v(c)\geq 1$ since $ C\cap A_{\leq 0} =
K$. Suppose that there exists an element $d\in C'$ with $\D (d)<\D
(c)$, we seek a contradiction. Let $q:= v(d)$. Then $q\geq 1$
since $ C\cap A_{\leq 0}=K$. Clearly, $v(c^q) = pq=v(d^p)$. Then,
by condition 4, $c^q= \l d^p +e$ for some nonzero scalar $\l$ and
an element $e\in C$ such that $v(e) < v(c^q)$. Using condition 2
we have the following strict inequalities:
\begin{eqnarray*}
 v(\d (c^q)) &=& v(c^q) +\D (c^q) = v (\l d^p) +\D (c) > v(\l d^p) +\D (d)= v(\l d^p) + \D (d^p) \\
              &=& v(\l d^p) + \D (\l d^p) = v(\d (\l d^p)),\\
               v(\d (c^q)) &=& v(c^q) +\D (c^q) = v (c^q) +\D (c) > v(e) +\D (e)= v(\d (e)).
\end{eqnarray*}
Now, $ v(\d (c^q)) >\max \{ v(\d (\l d^p)),  v(\d (e))\} \geq v(\d
(\l d^p)+ \d (e))=v(\d (\l d^p+e))= v(\d (c^q)), $
 a contradiction.

2. We keep the notation of the proof of statement 1. Let
\begin{equation}\label{mudet}
\mu := \max \{ h_k \, | \, k=0, 1, \ldots , mg^{-1}-1\}.
\end{equation}
For all natural numbers $j\geq \mu g^{-1}$, 
\begin{equation}\label{mudet1}
\dim_K(C_{\leq jg}) = j+1 -\nu \;\; {\rm where}\;\; \nu := |g\N
\backslash H|.
\end{equation}
Notice that the number  $g$ divides the drop $\D$ of the map $\d$.
For each natural number $j\geq \mu g^{-1}$, there is the short
exact sequence of vector spaces
$$ 0\ra K \ra C_{\leq jg}\stackrel{\d} {\ra}C_{\leq jg+\D} \ra
C_{\leq jg+\D}/ \d (C_{\leq jg })\ra 0,$$ and so, by
(\ref{mudet1}), 
\begin{equation}\label{mudet2}
\dim_K(C_{\leq jg+\D}/ \d (C_{\leq jg }))= \dim_K(C_{\leq jg+\D })
- \dim_K(C_{\leq jg}) +1= \frac{\D}{g}+1.
\end{equation}
Fix a subspace $U$ of $C$ such that $C= U\bigoplus \im (\d )$, and
so $U\simeq \coker (\d )$. Let $U'$ be a finite dimensional
subspace of $U$. Then $U'\subseteq C_{\leq jg+\D }$ for  some $j$,
hence $U'\bigoplus \d (C_{\leq jg})\subseteq C_{\leq jg+\D}$ and,
by (\ref{mudet2}), $\dim_K(U') \leq  \frac{\D}{g}+1$. This means
that $\dim_K( \coker (\d )) \leq \frac{\D}{g}+1$, i.e. the map
$\d$ is Fredholm since $\dim_K(\ker (\d )) = 1$, by condition 5.
Let 
\begin{equation}\label{gmud}
\g := \mu +\D .
\end{equation}
Then $C_{\leq \g } +\d (C_{\leq jg}) = C_{\leq jg+\D }$ for all
$j\geq j_0:= \mu g^{-1}$ since $\d$ is the map of constant drop
$\D$. The ascending chain of vector spaces $\{ W_j:= C_{\leq \g }
\cap \d (C_{\leq jg})\}_{j\in \N}$ of the finite dimensional
vector space $C_{\leq \g}$ stabilizers say at step $p$ and let
$W=W_p$. Choose a complementary subspace, say $V$, to $W$ in
$C_{\leq \g }$, i.e. $C_{\leq \g } = V\bigoplus W$. Then, for all
$j\geq \max \{ p , j_0\}$, $C_{\leq jg +\D } = V\bigoplus \d
(C_{\leq jg}).$ Since $C=\bigcup_{j\in \N} C_{\leq jg}$, we must
have $C= V\bigoplus \im (\d )$, and so $\dim_K(\coker (\d )) =
\dim_K(V) = \dim_K( C_{jg +\D } / \d (C_{\leq jg}) ) =
\frac{\D}{g} +1$, by (\ref{mudet2}). Then, $\ind (\d ) =
-\frac{\D}{g}$.

3. If $\D := \D(a)<0$ for some element $a\in C'$ then
 $v(\d (c))= v(c) +\D<v(c)$ for all elements
 $ c\in C'$. By conditions 3 and 5, the linear map $\d$ is locally nilpotent with kernel $K$.
 Since $g|\D$ and $\dim_K(\ker (\d )) =1$, we must have $\D =- g$
 and $H = \N g$; moreover, $\d ( C_{\leq jg}) = C_{\leq (j-1) g}$
 for all $j\in \N$ where $C_{\leq -g}:=0$. Then we can find a
 $K$-basis $\{ e_i\}_{i\in \N}$ for the algebra $C$ such that $\d
 (e_i) = e_{i-1}$ and $C_{\leq ig}= \bigoplus_{j=0}^i Ke_j$.

4. For all $a\in C\backslash \ker (\d \d')$, $\d'(a) \not\in \ker
(\d )$ and so $v(\d\d'(a)) = v( \d'(a)) +\D_{\d , v}= v(a)
+\D_{\d', v}+\D_{\d , v}$. Statement 4 follows also at once from
statement 2 and the additivity of the index.

 The proof of the
theorem is complete. $\Box $

$\noindent $

Let $D$ be a ring with an automorphism $\s $ and a central element
$a$.
 The {\bf generalized Weyl algebra} $A=D(\sigma, a)$  of degree 1
is the ring generated by $D$ and two indeterminates $X$ an $Y$
subject to the defining relations \cite{Bav-FA-1991},
\cite{Bav-AlgAnal-1992}:
$$
X\alpha=\sigma(\alpha)X \ {\rm and}\
Y\alpha=\sigma^{-1}(\alpha)Y,\; {\rm for \; all }\; \alpha \in D,
\ YX=a \ {\rm and}\ XY=\sigma(a).
 $$
The algebra $A={\bigoplus}_{n\in \mathbb{Z}}\, A_n$ is a
$\mathbb{Z}$-graded algebra where $A_n=Dv_n=v_nD$,
 $v_n=X^n\,\, (n>0), \,\,v_n=Y^{-n}\,\, (n<0), \,\,v_0=1.$

Let $K[H]$ be a polynomial algebra  in a variable $H$ over the
field $K$, $\s :H\ra H-1$ be the $K$-automorphism of the algebra
$K[H]$ and $a=H$. The first Weyl algebra $A_1=K<X, Y\;|\;YX-XY=1>$
is isomorphic to the generalized Weyl algebra
$$A_1\simeq K[H](\s , H),\; X\lra X,\; Y\lra Y,\; YX\lra H.$$
 We identify both these algebras via this isomorphism, that is
$A_1=K[H](\s , H)=\bigoplus_{i\in \Z}K[H]v_i$ and $H:=YX$.

$\noindent $ The Weyl algebra $A_1$ admits the following
$K$-algebra automorphism $\th$ and the anti-automorphism $\th'$
(i.e. $\th'(ab) = \th'(b) \th'(a)$ for all elements $a,b\in A_1$):
\begin{equation}\label{thA1}
\th : A_1\ra A_1, \;\; X\mapsto Y, \;\; Y\mapsto -X, \;\;
(H\mapsto -H+1);
\end{equation}
\begin{equation}\label{thA2}
\th' : A_1\ra A_1, \;\; X\mapsto Y, \;\; Y\mapsto X, \;\;
(H\mapsto H).
\end{equation}
They reverse the $\Z$-grading of the Weyl algebra $A_1$, i.e. $\th
(K[H]v_i) = K[H]v_{-i}$ and $\th' (K[H]v_i) = K[H]v_{-i}$
 for all elements $i\in \Z$.

$\noindent $

{\bf Proof of Theorem \ref{23Dec10}}. 3. Clearly, $K\subseteq \ker
(\d )$. Suppose that $\d (z) =0$ for some element $ z\in
C\backslash K$, we seek a contradiction. Recall that $C = C(z')$
for all elements $z'\in C\backslash K$ \cite{Dix} and $C(z')$ is a
finitely generated $K[z']$-module \cite{Ami}. In particular, $C$
is a finitely generated $K[z]$-module, and so every non-zero
element $c$ of $C$ is algebraic over $K[z]$, i.e. $f(c)=0$ for
some non-zero polynomial $f(t) \in K[z][t]$. We may assume that
its degree in $t$ is the least possible,  then $ f'(c)\neq 0$
where $f'=\frac{df}{dt}$, and  $0= \d ( f(c)) = f'( c) \d (c)$,
and so $\d (c)=0$. This means that $\d =0$, a contradiction.

1. We have to prove that the three statements are equivalent: $(a)
\Leftrightarrow (b) \Leftrightarrow (c)$. The implications
$(a)\Rightarrow (b)$ and $(c)\Rightarrow (a)$ are trivial.

$(b)\Rightarrow (c)$:  Suppose that the statement $(b)$ holds,
i.e. there exists an element $c\in C\backslash \ker (\d )$ such
that $\d^n (c)=0$ for some $n\geq 2$. Then $ 0\neq \d^m (c) \in
\ker (\d ) =K$ (statement 3) for some $m$ such that $ 1\leq m <n$.
Clearly, $t:= \d^{m-1} (c) / \d^m(c) \in C$ and $ \d (t) =1$.

$\noindent $

{\it Claim: $\d$ is a locally nilpotent derivation of the algebra
$C$}.

$\noindent $

Now, the implication $(b)\Rightarrow (c)$ follows from the Claim
and the well-known fact (which is a particular case of Theorem
\ref{2facts}.(1)): {\em Suppose that $A$ be a commutative algebra
over a field $K$ of characteristic zero, $\d \in \Der_K(A)$ and
$\d (t) =1$ for some element $t\in A$. Then $A$ is a polynomial
algebra $\ker (\d ) [ t]$ in $t$ with coefficients in the kernel
$\ker (\d )$ of $\d$}. In our situation, $\ker (\d) = K$
(statement 3), hence $C= K[t]$. Then $\d =\frac{d}{dt}$.

$\noindent $

{\bf Proof of the Claim}. Let $v$ be the $\Z$-graded degree on the
Weyl algebra $A_1$. In view of existence of grading reversing
automorphism $\th$  of the Weyl algebra $A_1$ see (\ref{thA1}), it
suffices to consider only two cases: either $t\in K[H]$ or $
v(t)>0$.

Suppose that  $t\in K[H]$.  Since $K[H]=C(H)$, and $t\in
C(H)\backslash K$, we have $C(H) = C(t)$. This means that $\d \in
\Der_K(K[H])$ and $ 1 = \d (t) = \frac{dt}{dH}\d (H)$. Therefore,
$\frac{dt}{dH} \d (H)\in K^*$, and so $t= \l H+\mu$ for some $\l
\in K^*$ and $\mu \in K$. Now, it is obvious that $K[H]= K[t]$ and
$\d = \frac{d}{dt}$.

Let $v(t) >0$ and $\D = \D_{\d , v}$. Notice that $\D (t) = v(1)
-v(t) = 0-v(t)= -v(t)<0$. Now, the Claim follows from the
following lemma.

\begin{lemma}\label{a23Dec10}
Suppose that $v(t)>0$. Then the conditions of Theorem
\ref{20Dec10} hold for $A=A_1$, $C$, $\d$, $\D = \D_{\d, v}$ with
$\D (t) = -v(t)<0$, and so $\d$ is a locally nilpotent derivation
of $C$.
\end{lemma}

{\it Proof}.  Condition 5 has been already established above.
Since $\d$ is a derivation of the commutative algebra $C$ over a
field of characteristic zero, conditions 1 and 2 of Theorem
\ref{20Dec10} hold by Lemma \ref{a17Dec10}. Notice that if
non-zero elements of the Weyl algebra $A_1$ commute then so do
their leading terms with respect to the $\Z$-grading of the Weyl
algebra $A_1$. The centralizer of all the homogeneous elements of
the Weyl algebra $A_1$ are found in Proposition 3.1,
\cite{BavDP5}. This description makes conditions 3 and 4  obvious.
The proof of Lemma \ref{a23Dec10} is complete. $\Box$

2. We have to prove that the four statements are equivalent: $
(a)\Leftrightarrow (b)\Leftrightarrow (c)\Leftrightarrow (d)$. The
implications $(a)\Rightarrow (b)$, $(c)\Rightarrow (a)$,
$(c)\Rightarrow (d)$ and $(d)\Rightarrow (a)$ are trivial. It
remains to show that the implication $(b)\Rightarrow (c)$ holds.

$(b)\Rightarrow (c)$ Suppose that $\d (c) = \l c$ for some $0\neq
c\in C$ and $\l \in K^*$. In view of existence of the grading
reversing automorphism $\th$  of the Weyl algebra $A_1$ see
(\ref{thA1}), it suffices to consider only two cases: either $c
\in K[H]$ or $v(c)>0$.

Suppose that  $c\in K[H]$. Since $\ker (\d ) =K$ (statement 3), we
see that $c\in C\backslash K$, and so $C= C(c) = C(H) = K[H]$. The
five conditions of Theorem \ref{20Dec10} are satisfied for $A= C =
K[H]$, $v = \deg_H $, $\d$ and $\D = \D_{\d , v}$ and so, for all
polynomials $p\in K[H]\backslash K$, $\deg_H(\d (p)) = \deg_H(p)$
since $\deg_H(\d (c)) = \deg_H(\l c)  =\deg_H(c)$. Then $\d = (\mu
H +\nu ) \frac{d}{dH}$ for some scalars $\mu \in K^*$ and $\nu \in
K$. Let $H':= \mu H +\nu$. Then $K[H] = K[H']$ and $\d =
H'\frac{dH'}{dH}\frac{d}{dH'}=\mu H'\frac{d}{dH'}$. Now, the
statement $(c).(i)$  is obvious as $\d (H'^i) = i\mu H'^i$ for all
$i\in \N$.

Let $v(c)>0$ and $\D = \D_{\d , v}$. Notice that $\D (c) = v(\l
c)-v(c) = v(c) - v(c) =0$.

\begin{lemma}\label{a25Dec10}
Suppose that $v(c)>0$. Then the conditions of Theorem
\ref{20Dec10} hold for $A=A_1$, $C$, $\d$, $\D = \D_{\d , v}$ with
$\D (c)=0$, and so the derivation $\d$ respects the filtration
$C=\bigcup_{i\in \N} C_{\leq i}$ where $C_{\leq i}:= \{ c\in C \,
| \, v(c) \leq i\}$, i.e. $\d (C_{\leq i}) \subseteq C_{\leq i}$
for all $i\in \N$.
\end{lemma}

{\it Proof}. Repeat word for word the proof of Lemma
\ref{a23Dec10} $\Box $

$\noindent $

By Lemma \ref{a25Dec10}, $C_0 = K$ and $\supp (C) \subseteq \N$.
By Lemma \ref{a25Dec10} and the fact that $\dim_K(C_{\leq i} /
C_{\leq i-1}) \leq 1$ for all $i\in \N$, the algebra $C$ is the
direct sum $C=\bigoplus_{\l \in \Ev (\d )}C^\l$ where $C^\l :=
\bigcup_{n\geq 1} \ker (  ( \d - \l )^n)$ where $\Ev (\d )$ is the
set of eigenvalues of the derivation $\d \in \Der_K(C)$ in the
field $K$.

$\noindent $

{\it Claim: $C^\l = \ker (\d - \l )$ for all $\l \in \Ev (\d )$.}

$\noindent $

To prove the claim, for each nonzero element $c_\l \in C^\l$,  we
introduce the {\em nilpotency degree} $ d_\l (c_\l )$  of $c_\l$
 by the rule
$$ d_\l (c_\l ) := \min \{ n\in \N \, | \, (\d -\l)^{n+1}(c_\l )=0\}.$$
 For all $\l , \mu \in \Ev (\d )$, $0\neq c_\l
 \in C^\l$ and $ 0\neq c_\mu \in C^\mu$,
\begin{equation}\label{dccm}
d_{\l +\mu } (c_\l c_\mu ) = d_\l ( c_\l ) + d_\mu (c_\mu ).
\end{equation}
In more detail, let $n= d_\l (c_\l )$ and $m = d_\mu (c_\mu )$. It
follows from the equality: for all elements $a,b\in C$,
\begin{equation}\label{dccm1}
(\d -\l - \mu )^i (ab) = \sum_{i=0}^n {n\choose i} (\d - \l )^i
(a) (\d -\mu )^{n-i}(b)
\end{equation}
that $(\d -\l - \mu ) ^{n+m+1} (c_\l c_\mu ) =0$ and $(\d -\l -
\mu )^{n+m} (c_\l c_\mu ) = {n+m\choose n}(\d - \l )^n (c_\l )
\cdot (\d -\mu )^m (c_\mu ) \neq 0$. This proves (\ref{dccm}).

Recall the Theorem of Amitsur \cite{Ami}: {\em for any element
$c\in C\backslash K$, the algebra $C$ is a finitely generated free
$K[c]$-module}. In particular, the Gelfand-Kirillov dimension of
the algebra $C$ is 1. Suppose that $C^\l \neq \ker (\d -\l )$ for
some $\l \in \Ev (\d )$, we seek a contradiction. By (\ref{dccm}),
we may assume that $\l \neq 0$. Fix  nonzero elements, say $s,t\in
C^\l$, such that $s\in \ker (\d -\l )$, $t\not\in \ker (\d - \l )$
but $(\d - \l )^2 (t)=0$. Then the elements $s $ and $t$ are
algebraically independent: if $ \sum \l_{ij} s^it^j=0$ for some
scalars $\l_{ij}\in K$ not all of which are equal to zero then we
may assume that, for all $i$ and $j$, $i+j=m=const$ since
$s^kt^l\in C^{(k+l)\l }$ for all $k$ and $l$; let $p:= \max \{ j
\, | \, \l_{m-j, j}\neq 0\}$; then, by (\ref{dccm1}),
\begin{eqnarray*}
0&= & (\d -m\l )^p(\sum_{j=0}^p\l_{m-j, j} s^{m-j}t^j) =
 \sum_{j=0}^p\l_{m-j, j}(\d -(m-j)\l -  j\l )^p (s^{m-j}t^j) \\
 &= & \sum_{j=0}^p\l_{m-j, j}\sum_{i=0}^p{p\choose i} (\d -
 (m-j)\l ) ^{p-i}(s^{m-j}) (\d - j\l )^i (t^j)= \sum_{j=0}^p\l_{m-j,
 j}s^{m-j}(\d -j\l )^p(t^j)\\
 &=& \l_{m-p, p}s^{m-p}(\d - p \l)^p (t^p) = \l_{m-p, p}s^{m-p}\cdot
 p!\cdot
 ((\d -\l )(t))^p\neq 0,
\end{eqnarray*}
a contradiction.  Since the elements $s$ and $t$ of the
commutative  algebra $C$ are algebraically independent its
Gelfand-Kirillov is at least 2, a contradiction. The proof of the
Claim is complete.

Let $\l ,\mu \in \Ev (\d ) \backslash \{ 0\}$, $0\neq c_\l \in
\ker (\d - \l )$ and $0\neq c_\mu \in \ker ( \d - \mu )$. By Lemma
\ref{a25Dec10}, $C_0 = K$, $\ker (\d ) =K$,  and so $p:= v(c_\l
)\geq 1$, $q:= v(c_\mu )\geq 1$, and there exists a scalar $\nu\in
K^*$ such that $c_\mu^p= \nu c_\l ^q+\cdots $ where the three dots
denote smaller terms (i.e. of smaller $\Z$-graded degree than
$v(\nu c_\l^q)=pq$). Applying the derivation $\d$ to this equality
and using Lemma \ref{a25Dec10}, we see that $p\mu c_\mu^p = \nu
q\l c_\l^q+\cdots$. Comparing the leading terms of the two
equalities we obtain the equality $$\mu = \frac{\l}{p}q=
\frac{\l}{v(c_\l )}v(c_\mu ).$$ Therefore, $\Ev (\d )
=\frac{\l}{v(c_\l )}\supp (C)$, by the Claim, and the case $(ii)$
follows. The proof of Theorem \ref{23Dec10} is complete. $\Box $

In general, it is not true that in cases $2(ii), (iii)$ of Theorem
\ref{23Dec10}, the algebra $C$ is isomorphic to a polynomial
algebra or is a non-singular algebra.

$\noindent $

{\it Example 4}. Let $v = H(H-1)^{-1} (H-2) X\in \Frac (A_1)$.
Then $v\not\in A_1$ but $v^i\in A_1$ for all $i\geq 2$. By
Proposition
 3.1, \cite{BavDP5}, $C:= C(v^2) = K[v]\cap A_1= K\bigoplus
 \bigoplus_{i\geq 2} Kv^i$. Therefore, $\supp (C) = \{ 0, 2,3,
 \ldots \}\neq \N$. For all $i\geq 0$, $[H, v^i]=iv^i$. Therefore,
 the restriction $\d$ of the inner derivation $\ad (H)$ to the
 $\ad (H)$-invariant algebra $C$ yields a semi-simple derivation
 with the set of eigenvalues $\{ 0,2,3,\ldots \}$, and the algebra
 $C$ is not isomorphic to a polynomial algebra. The algebra $C$ is
 isomorphic to the algebra $K[s,t]/(s^2-t^3)$ of regular functions
 on the cusp $s^2= t^3$. In particular, the algebra $C$ is a
 singular one.

$\noindent $

For each element $a\in A_1$, the union $N(a) := N(a, A_1) :=
\bigcup_{i\geq 0} N(a, A_1, i)$, where $N(a,A_1, i):= \ker ( \ad (
a)^{i+1})$, is a filtered algebra ($N(a,A_1,i) N(a,A_1, j)
\subseteq N(a, A_1,i+j)$ for all $i,j\geq 0$). By the very
definition, the algebra $N(a)$ is the largest subalgebra of the
Weyl algebra $A_1$ on which the inner derivation $\ad (a)$ acts
locally nilpotently. Little is known about these algebras in the
case when $N(a) \neq C(a)$. In particular, it is not known of
whether these algebras are finitely generated or Noetherian.
Though, a positive answer to Dixmier's Fourth Problem \cite{Dix},
which is still open,  would imply that the algebras $N(a)$ are
finitely generated and Noetherian. In case of {\em homogeneous}
elements of the Weyl algebra $A_1$, a positive answer to Dixmier's
Fourth Problem was given in \cite{BavDP5}. In particular, for all
homogeneous elements $a$ of $A_1$, the algebra $N(a)$ is a
finitely generated and Noetherian.

\begin{proposition}\label{a18Dec10}
If $[y,x]=1$ for some elements $x,y\in A_1$ then $N(x, A_1) = N(y,
A_1) = K\langle x, y \rangle$.
\end{proposition}

{\it Proof}.  In view of existence of the $K$-algebra automorphism
of the Weyl algebra $A_1':=K\langle x, y \rangle \ra A_1'$,
$x\mapsto y$, $y\mapsto -x$, it suffices to show that the algebra
$N:= N( y , A_1)$ is equal to $A_1'$. The inner derivation $\d =
\ad (y)$ is a locally nilpotent derivation of the algebra $N$ with
$\d (x) = 1$. By Theorem \ref{2facts}.(1)  and Theorem
\ref{GGV-Cen}, $N= \ker (\d ) [ x; \ad (x)]= C(y) [ x; \ad (x) ] =
K[y][x; \ad (x) ] = A_1'$. $\Box $

$\noindent $

For each element $a\in A_1\backslash K$, let $\Ev (a)$ be the set
of eigenvalues in the field $K$ of the inner derivation $\ad (a)$
of the Weyl algebra $A_1$. Then $\Ev (a)$ is an additive submonoid
of $(K, +)$ and 
\begin{equation}\label{DaEv}
D(a) := \bigoplus_{ \l \in \Ev (a)} D( a, \l), \;\; {\rm
where}\;\; D(a, \l ) := \ker_{A_1} (\ad (a) - \l ),
\end{equation}
is an $\Ev (a )$-graded subalgebra of the Weyl algebra $A_1$, i.e.
$$ D(a, \l ) D(a, \mu ) \subseteq D(a, \l +\mu ) \;\; {\rm for \;
all }\;\; \l , \mu \in \Ev (a).$$ Little is known about the
algebras $D(a)$ in general.

$\noindent $

The next corollary is a first  step in the proof of Theorem
\ref{22Dec10}.

\begin{corollary}\label{a19Dec10}
Suppose that  $[y,x]=1$ for some elements $x,y\in A_1$. Let
$A_1':= K\langle x,y \rangle$, $\d_x := \ad (x)$, $\d_y := \ad
(y)$,  and $h:= yx$. Then $D(h) = A_1'$ iff $\d_x$ is a locally
nilpotent derivation  of the algebra $D(h)$ iff $\d_y$ is a
locally nilpotent derivation of the algebra $D(h)$ iff $\d_x\d_y$
is a locally nilpotent map in $D(h)$.
\end{corollary}

{\it Proof}. If $D(h) = A_1'$ then $\d_x$ is a locally nilpotent
derivation of the algebra $D(h)$. Conversely, if $\d_x$ is a
locally nilpotent derivation of the algebra $D(h)$ then $D(h)
\subseteq N(x, A_1) = A_1'$, by Proposition \ref{a18Dec10}.
Therefore, $D(h) = A_1'$ since the inclusion $A_1'\subseteq D(h)$
is obvious. By symmetry, the second `iff' is also true. The
derivations $\d_x$ and $\d_y$ commute. So, if $D(h) = A_1'$ then
the map $\d := \d_x\d_y$ is a locally nilpotent map on $D(h)=
A_1'$ since  $\d_x$ and $\d_y$ are commuting locally nilpotent
derivations of the Weyl algebra $A_1'$. Conversely, if $\d$ is a
locally nilpotent map on $D(h)$ then for any element $a\in D(h)$,
$0=\d^n (a) = \d_x^n \d_y^n (a) $ for some natural number $n\geq
1$. Thus $\d_y^n (a) \in N (x, A_1) = A_1' = N(y, A_1)$ (by
Proposition \ref{a18Dec10}), and so $a\in N(y, A_1) = A_1'$. This
means that $D(h) = A_1'$. $\Box $


\begin{proposition}\label{A18Dec10}
Suppose that $[y,x]=1$ for some elements $x,y\in A_1$. Let $h:=
yx$,  $\d_x:= \ad (x)$, $\d_y := \ad (y)$, $\d := \d_x\d_y$ and $
A_1':= K\langle x,y\rangle$. Then
\begin{enumerate}
\item $N(\d , A_1) = A_1'$. \item $\ker_{A_1}(\d ) = K[x]+K[y]$.
\item $\ker_{C(h)}(\d ) = K$. \item $C(h) \cap A_1'= K[h]$. \item
$C(h) = K[h]$ iff the map $\d$ acts locally nilpotently on $C(h)$.
\item $K[h]=\bigoplus_{i\geq 0} K \frac{y^ix^i}{(i!)^2}=
\bigoplus_{i\geq 0} K \frac{x^iy^i}{(i!)^2}$,  $\d
((-1)^i\frac{y^ix^i}{(i!)^2})=(-1)^{i-1}\frac{y^{i-1}x^{i-1}}{((i-1)!)^2}$
and
 $\d
((-1)^i\frac{x^iy^i}{(i!)^2})=(-1)^{i-1}\frac{x^{i-1}y^{i-1}}{((i-1)!)^2}$
for all $i\geq 1$ where $\frac{y^ix^i}{(i!)^2}=\frac{h(h+1)\cdots
(h+i-1)}{(i!)^2}$ and
$\frac{x^iy^i}{(i!)^2}=\frac{(h-1)(h-2)\cdots (h-i)}{(i!)^2}$.
\end{enumerate}
\end{proposition}

{\it Proof}. 1. Since the derivations $\d_x$ and $\d_y$ commute
and $N(x, A_1) = N( y , A_1) = A_1'$ (by Proposition
\ref{a18Dec10}), the inclusion $A_1'\subseteq N(\d , A_1)$
follows. Let $a\in N( \d , A_1)$. Then $0=\d^n(a) = \d_x^n\d_y^n
(a)$ for some $n\geq 1$, hence $\d_y^n (a) \in N(\d_x, A_1) =
A_1'= N(y, A_1)$ (by Proposition \ref{a18Dec10}), and so $a\in
N(\d_y, A_1) = A_1'$, i.e.  $N(\d , A_1) = A_1'$.

2. Clearly, $K[x]+K[y]\subseteq \ker _{A_1} (\d ) \subseteq N(\d ,
A_1) = A_1'$, by statement 1. Therefore, $\ker_{A_1}(\d ) =
\ker_{A_1'}(\d )$. Since $\d (x^iy^j) = -ij x^{i-1}y^{j-1}$ for
all $i,j\geq 1$, we have the opposite inclusion $\ker_{A_1} (\d )
\subseteq K[x]+K[y]$.

3. $\ker_{C(h)}(\d ) =C(h) \cap \ker_{A_1}(\d )= C(h) \cap
(K[x]+K[y])=K$, by statement 2 and the fact that $[h,x^i]=ix^i$
and $ [ h, y^i]=-iy^i$ for all $i\geq 0$.

4. The Weyl algebra $A_1'= \bigoplus_{i\in \Z}K[h]v_i'$ is a
$\Z$-graded algebra where
$$ v_i:= \begin{cases}
x^i& \text{if }i\geq 0,\\
y^{-i}& \text{if }i<0.\\
\end{cases}$$
and $ [ h,u]=iu$ for all $u\in K[h]v_i'$. Therefore, $C(h) \cap
A_1'= K[h]$.

5. $(\Rightarrow )$ This follows from $\d_x^2(h) = 0$, $\d_y^2(h)
= 0$ and $\d_x\d_y = \d_y\d_x$ (and so $\d_x^{n+1} (h^n) =
\d_y^{n+1} (h^n)=0$ for all $n\geq 1$).

$(\Leftarrow )$ If $\d$ acts locally nilpotently on $C(h)$ then
$C(h) \subseteq N(\d , A_1) = A_1'$, by statement 1, and so $C(h)
= C(h) \cap A_1' = K[h]$, by statement 4.

6. Statement 6 follows at once from the following two facts:
$\deg_h(y^ix^i) = \deg_h(x^iy^i) = i$ for all $i\geq 0$ and $\d
(y^ix^i ) = -i^2y^{i-1}x^{i-1}$ and $\d (x^i y^i ) = -i^2
x^{i-1}y^{i-1}$ for all $i\geq 1$.  $\Box $

$\noindent $

{\bf The degree function $v_{\rho , \eta}$ on the Weyl algebra
$A_1$}.  The elements $\{ Y^iX^i\, | \, (i,j) \in \N^2\}$ is a
$K$-basis of the Weyl algebra $A_1$. Any pair $(\rho , \eta )$ of
positive integers determines the degree function $v_{\rho , \eta
}$ on the Weyl algebra $A_1$ by the rule 
\begin{equation}\label{vrhe}
v_{\rho , \eta } (\sum_{i,j\in \N} a_{ij}Y^iX^j) := \max \{ \rho i
+\eta j \, | \, a_{ij}\in K^*\}.
\end{equation}
Since $\im (v_{\rho , \eta }) \subseteq  \N$, the negative terms
of the $v_{\rho , \eta }$-filtration are all equal to zero, i.e.
$A_1=\bigcup_{i\in \N} A_{1, \leq i } (\rho , \eta )$ is a
positively filtered algebra where $A_{1, \leq i } (\rho , \eta ):=
\{ a\in A_1\, | \, v_{\rho , \eta }(a) \leq i\}$, and $A_{1, \leq
0} (\rho , \eta ) = K$. It follows from the relation $[Y,X]=1$
that, for all elements $ a,b\in A_1$, 
\begin{equation}\label{vrab}
v_{\rho , \eta }([a,b]) \leq v_{\rho , \eta }(a) + v_{\rho , \eta
}(b) - \rho - \eta ,
\end{equation}
and so the associated graded algebra ${\rm gr}_{ v_{\rho , \eta }}
(A_1)$ is a polynomial algebra in two variables which are the
images of the elements $X$ and $Y$ in ${\rm gr}_{ v_{\rho , \eta
}} (A_1)$.

For each non-zero element $a= \sum_{i,j\in \N} a_{ij}Y^iX^j$ where
$a_{ij}\in K$, define its {\em Newton polygon} $\NP (a)$ as the
the convex hull of the {\em support} $\supp (a) := \{ (i,j)\in
\N^2\, | \, a_{ij}\neq 0\}$ of the element $a$. The pair $(\rho,
\eta )$ is called $a$-{\em generic} if 
\begin{equation}\label{rege}
\# \{ (i,j) \in \supp (a) \, | \rho i +\eta j = v_{\rho , \eta}
(a)\} =1.
\end{equation}
 The set of $a$-generic pairs is a non-empty set (moreover, all
 but finitely many pairs $(\rho , \eta)$ are $a$-generic, eg,
 pairs such that the lines $\rho y +\eta x =1$ are not parallel to
 the edges of the Newton polygon of the element $a$).

$\noindent $

{\bf Proof of Theorem \ref{X1Jan11}}. Let $v=v_{\rho , \eta}$
where $\rho$ and $\eta$ are positive integers, $A_1'':= \Frac
(A_1') \cap A_1$. The algebras $A_1$, $A_1'$ and $A_1''$ are
invariant under the action of the inner derivation $\d_x := \ad
(x)$.

Clearly, $A_1'\subseteq A_1''$. Suppose that $A_1'\neq A_1''$, we
seek a contradiction. Fix an element $a\in A_1''\backslash A_1'$
with the least possible $v(a)$. There exist nonzero elements
$p,q\in A_1'$ such that $pa=q$. Let $ \d_x(\cdot ) = (\cdot )'$.
Then, by (\ref{vrab}),  $v(p')<v(p)$, $v(q') < v(q)$, and
$q'=p'a+pa'$, and so
$$ v(p)+v(a')= v(pa')=v(q'-p'a)\leq \max \{ v(q'), v(p')+v(a)\} <
\max \{ v(q), v(p)+v(a)\} = v(p)+v(a).$$ Therefore, $v(a') < v(
a)$ and $ a'\in A_1''$. By the minimality of $v(a)$, $a'\in A_1'=
N(\d_x,  A_1)$ (Proposition \ref{a18Dec10}), hence $a\in N(\d_x,
A_1) = A_1'$, a contradiction. $\Box $

$\noindent $

Till the end of this section, we assume that $[y,x]=1$ for some
 elements $x$ and $y$ of the Weyl algebra $A_1$ and $h:= yx$. The
 results below are some of the key steps in the proof of Theorems
 \ref{18Dec10} and \ref{22Dec10}.

Let $d := \d_y (\cdot ) x: A_1\ra A_1$, $a\mapsto [ y,a]x$. Then,
for all elements $a,b\in A_1$,  
\begin{equation}\label{dab}
d(ab) = d(a) b + ad(b) +d(a) h^{-1}y[b,x].
\end{equation}
In more detail, using the fact that that $\d_y$ is a derivation of
the Weyl algebra $A_1$, we obtain the equality
\begin{eqnarray*}
 d(ab) &=& \d_y(a)xy\cdot (xy)^{-1}  bx+ ad (b) = d (a) y \frac{1}{h-1} bx + ad(b)\\
 &=& d(a)\frac{1}{h}ybx   + ad(b)=d(a) b + ad(b) +d(a) ( h^{-1}y
 bx-b)\\
 &=&d(a) b + ad(b) +d(a) h^{-1}y [b,x].
\end{eqnarray*}
Similarly, let $d' := \d_x (\cdot ) y: A_1\ra A_1$, $a\mapsto [
x,a]y$. Then, for all elements $a,b\in A_1$,  
\begin{equation}\label{sdab}
d'(ab) = d'(a) b + ad'(b) +d'(a) (h-1)^{-1}x[b,y].
\end{equation}
In more detail, using the fact  that $\d_x$ is a derivation of the
Weyl algebra $A_1$, we obtain the equality
\begin{eqnarray*}
 d'(ab) &=& \d_x(a)yx\cdot (yx)^{-1}  by+ ad' (b) = d' (a) x h^{-1} by + ad'(b)\\
 &=& d'(a)\frac{1}{h-1}xby   + ad'(b)=d'(a) b + ad'(b) +d'(a) (
 (h-1)^{-1}x
 by-b)\\
 &=&d'(a) b + ad'(b) +d'(a) (h-1)^{-1}x [b,y].
\end{eqnarray*}

Let $d\in \{ \d_y(\cdot )x, \d_x (\cdot ) y, x\d_y(\cdot ) , y
\d_x (\cdot ) \}$ . Then $d(C(h)) \subseteq C(h)$ since
$$d(C(h)) \subseteq D(h,1)D(h,0)D(h, -1)+D(h,-1)D(h,0)D(h, 1)\subseteq D(h, 1+0-1)=
D(h, 0)= C(h).$$
\begin{lemma}\label{a21Dec10}
Let $d:= \d_y (\cdot ) x$, $v= v_{\rho , \eta}$ where $ \rho $ and
$\eta $ are positive integers. Then $v(d(a^n)) = v (a^{n-1} d(a))$
for all elements $ a\in C(h) \backslash K$ and $n\geq 1$.
\end{lemma}

{\it Proof}. Notice that $\ker_{C(h)} (d) = C(h) \cap C(y) = K$
since the elements $h$ and $y$ do not commute. Therefore, if $a\in
C(h) \backslash K$ then $d(a^n) \neq 0$ for all $n\geq 1$. Let us
prove by induction on $n$ that 
\begin{equation}\label{dana}
d(a^n) = na^{n-1} d(a) +\cdots \;\; {\rm for \; all}\;\; a\in C(h)
\backslash K,
\end{equation}
where the three dots denote smaller terms, i.e. $v(\cdots )
<v(na^{n-1} d(a))$. The claim is trivially true for $n=1$. So, let
$n>1$. By (\ref{dab}),
$$ d(a^n) = d(a\cdot a^{n-1}) = d(a) a^{n-1} +ad(a^{n-1}) + d(a)h^{-1}y[a^{n-1},x].$$
We can extend the degree function $v$ from the Weyl algebra $A_1$
to its quotient ring $\Frac (A_1)$, the, so-called, {\em Weyl skew
field} of fractions of $A_1$, by the rule $v(s^{-1} a) = v(a)
-v(s)$. Then, by (\ref{vrab}),
\begin{equation}\label{dana1}
v(d(a) h^{-1} y[a^{n-1}, x]) \leq v( d(a)) - v(h) +v(y)
+v(a^{n-1})  +v(x) - \rho -\eta = v(a^{n-1} d(a))-\rho -\eta .
\end{equation}
By induction on $n$, $a\cdot  d(a^{n-1}) = (n-1) a\cdot a^{n-2}
d(a) +\cdots =(n-1)a^{n-1} d(a)+\cdots $. Recall that the algebra
$C(h)$ is commutative. In particular, $d(a)a^{n-1} = a^{n-1}d(a)$.
Then
\begin{eqnarray*}
d(a^n) &= & a^{n-1} d(a) +((n-1) a^{n-1} d(a) +\cdots )
+d(a) h^{-1} y[a^{n-1}, x]  \\
& =&(na^{n-1} d(a) +\cdots )+d(a)h^{-1} y[a^{n-1}, x]   = na^{n-1}
d(a) +\cdots \;\;\; ({\rm by} \;\; (\ref{dana1})). \;\; \Box
\end{eqnarray*}

\begin{theorem}\label{A20Dec10}
Let $v=v_{\rho , \eta}$ where $(\rho , \eta )$ is an $h$-generic
pair and let $d\in \{ \d_y(\cdot )x, \d_x (\cdot ) y, x\d_y(\cdot
) , y \d_x (\cdot ) \}$ where $\rho$ and $ \eta$ are positive
integers. Then $\D _{d, v} (c) =0$ for all elements $c\in C(h)
\backslash K$. Moreover, $A=A_1$, $C=C(h)$, $v$, $\d := d$ and $\D
= \D_{\d, v}$ satisfy the five conditions of Theorem
\ref{20Dec10}.
\end{theorem}

{\it Proof}. In view of existence of the $K$-algebra isomorphism
of the Weyl algebra  $A_1':= K\langle x,y\rangle \ra A_1'$,
$x\mapsto y$, $y\mapsto -x$  ($h=yx\mapsto -xy =-h+1$), and the
$K$-algebra anti-isomorphism $A_1' \ra A_1'$, $x\mapsto y$,
$y\mapsto x$  ($h=yx\mapsto yx =h$), and the fact that $C(h) =
C(-h+1)$, it suffices to prove the theorem only for the map  $d=
\d_y(\cdot ) x$. Notice that the fifth condition of Theorem
\ref{20Dec10} holds: $\ker_{C(h)}(d) = \ker_{C(h)} (\d_y) = K$, by
Proposition \ref{A18Dec10}.(3) (since $K\subseteq \ker_{C(h)}
(\d_y) \subseteq \ker_{C(h)} (\d_x\d_y ) = K$); and $d(h) =
[y,h]x= yx =h$. Therefore, $\D (h) =0$. Now, to finish the proof
of the theorem it suffices to show that $A_1$, $C(h)$, $v$, $d$
and  $\D = \D_{d, v}$ satisfy the first four conditions of Theorem
\ref{20Dec10}.

 The first condition of Theorem \ref{20Dec10} follows from
 (\ref{vrab}) and (\ref{dab}): for all elements $a, b\in C(h)$,
$$v(d(a) h^{-1} y[b,x] ) \leq  v( d(a)) -v(h) +v(y) +v(b)+ v(x) -\rho - \eta =
v( d(a) b) -\rho -\eta <v( d(a) b).$$ Then, by (\ref{dab}),
\begin{eqnarray*}
 \D (ab) &=& v( d(ab))-v(ab)\leq \max \{ v(d(a)b+d(a)h^{-1} y[b,x]), v( ad(b))\}-v(a) - v(b)  \\
 &=& \max \{  v(d(a)b)-v(a) - v(b), v( ad(b))-v(a) - v(b)= \max \{ \D (a) ,
 \D (b)\}.
\end{eqnarray*}
The second condition of Theorem \ref{20Dec10} follows from Lemma
\ref{a21Dec10}: for all $a\in C(h) \backslash K$ and $n\geq 1$,
\begin{eqnarray*}
 \D (a^n) &=& v( d(a^n))-v(a^n)= v(a^{n-1} d(a)) - nv(a) = (n-1)
 v(a) +v(d(a))-nv(a)\\
 &=& v(d(a)) - v(a) = \D (a);
\end{eqnarray*}
for all $a\in K$ and $n\geq 1$, $\D (a^n ) =-\infty = \D (a)$. The
third condition of Theorem \ref{20Dec10} is obvious.

Let $(i,j) $ be the only pair satisfying (\ref{rege}) for the
element $h$. Then $h = \l Y^iX^j+\cdots $ for some $\l \in K^*$
where the three dots mean smaller terms with respect to the
$v$-filtration on the Weyl algebra $A_1$. The element $l(h) := \l
Y^iX^j$ is called the {\em leading term} of the element $h$ with
respect to the $v$-filtration. Since
$$ Y^i X^j = \begin{cases}
(H+i-1) (H+i-2) \cdots (H+1) H X^{j-i} & \text{if }i\leq j,\\
(H+i-1) (H+i-2) \cdots (H+i-j) Y^{i-j}& \text{if }i>j,\\
\end{cases}$$
the leading element $l(h)$ is a {\em homogeneous} element of the
$\Z$-graded algebra $A_1$. The centralizers  of all homogeneous
elements of $A_1$ are found explicitly Proposition 3.1,
\cite{BavDP5}. It follows easily from this result and the obvious
fact that $[h,h']=0$ implies $[l(h), l(h')]=0$ that the fourth
condition of Theorem \ref{20Dec10} holds. $\Box $


\begin{corollary}\label{aA20Dec10}
Let the degree function $v$ on the Weyl algebra $A_1$ be as in
Theorem \ref{A20Dec10}, $\d := \d_x\d_y$ and $\D = \D_{\d , v}$.
Then $\D (c) = -v(h) <0$ for all elements $ c\in C(h)\backslash
K$.
\end{corollary}

{\it Proof}. Let $d:= \d_y(\cdot ) x$ and $d':= \d_x(\cdot ) y$.
Then 
\begin{equation}\label{ddp}
dd'= \d_y(\d_x(\cdot )y) x= \d_y\d_x(\cdot ) yx = \d (\cdot ) h.
\end{equation}
Recall that $\ker_{C(h)}(\d) = K$ (Proposition
\ref{A18Dec10}.(3)), $\ker_{C(h)}(d) =C(h) \cap C(y)=K$ since
$hy\neq yh$, and $\ker_{C(h)}(d') =C(h) \cap C(x)=K$ since $hx\neq
xh$. Therefore, if $\d (c) \neq 0$ for some element $c\in C(h)$
then $ d' (c) \not\in K$. By Theorem \ref{20Dec10}.(4) and Theorem
\ref{A20Dec10}, for all elements $c\in C(h) \backslash K$,
$v(dd'(c)) = v(c) +\D_{d, v} +\D_{d', v}=v(c) +0+0=v(c)$, and so
$\D_{dd', v}(c)=0$. On the other hand, for all elements $c\in C(h)
\backslash K$,
$$ 0=\D_{dd', v}(c) = \D_{\d (\cdot )h, v}(c) = v(\d (c) h) -
v(c) = \D (c) +v(h),$$ and so $\D (c) = -v(h) <0$. $\Box $

$\noindent $

{\bf Proof of Theorem \ref{18Dec10}}.  By Corollary
\ref{aA20Dec10}, the map $\d$ acts locally nilpotently on the
algebra $C(h)$ since,  for all elements $c\in C(h) \backslash K$,
$ v(\d (c)) = v(c) +\D (c) = v(c) -v(h)<v(c)$ and $v(c')\geq 0$
for all elements $ c'\in C$ and $\d (K)=0$. Then, $C(h) = K[h]$,
by Proposition \ref{A18Dec10}.(5). $\Box$.


\section{Proof of Theorem \ref{22Dec10}}\label{PR22DEC}

The aim of this section is to give the proof of Theorem
\ref{22Dec10}. In this section, we assume that $[y,x]=1$ for some
elements $x$ and $y$ of the Weyl algebra $A_1$ and $h:= yx$. Let
$A_1':= K\langle x,y\rangle$, $\d_x:= \ad (x)$, $\d_y := \ad (y)$
and $\d := \d_x\d_y= \d_y\d_x$.

{\bf Proof of Theorem \ref{22Dec10}}.  Notice that the Weyl
algebra  $A_1':= K\langle x,y \rangle=\bigoplus_{i\in \Z}K[h]v_i'$
is a subalgebra of $D(h)$ since, for all elements $u\in K[h]v_i'$,
$[h,u]=iu$; and the set of integers $\Z$ is a subset of the set
$\Ev (h)$ of eigenvalues of the inner derivation $ \d := \ad (h)$
 of the Weyl algebra $A_1$. By the Theorem of Joseph
  (\cite{joseph-AJM}, Note added in proof), $\Ev (h) = \Z \th$
for some nonzero  scalar $\th $ necessarily of the form $\pm
\frac{1}{n}$ since $1\in \Z \th$. Without loss of generality we
may assume that $\th = \frac{1}{n}$. Then
$$D(h) =
\bigoplus_{i\in \Z} D_{i\th}\;\; {\rm where}\;\; D_{i\th }:= \{
a\in A_1\, | \, [ h, a]=i\th a\}.$$ For all $\l , \mu \in \Ev (
h)$, $D_\l D_\mu \subseteq D_{\l +\mu }$. By Theorem
\ref{18Dec10}, $D_0= C(h) = K[h]$.

$\noindent $

{\em Claim 1: For each $\l \in \Ev (h)$, there exists an element
$u_\l \in D_\l$ such that $D_\l = K[h]u_\l = u_\l K[h]$.}

$\noindent $

If $\l =0$ then take $u_0=1$, by Theorem \ref{18Dec10}. If $\l
\neq 0$ then take any  nonzero element, say $w_\l$, of $D_\l$ then
the map $\cdot w_{-\l} : D_\l \ra D_0=K[h]$, $ d\mapsto dw_{-\l}$,
is a monomorphism of left $K[h]$-modules. Therefore, $D_\l=
K[h]u_\l$ for some element $u_\l \in D_\l$. Similarly, the map $
w_{-\l} \cdot: D_\l \ra D_0=K[h]$, $ d\mapsto w_{-\l} d$,  is a
monomorphism of right  $K[h]$-modules. Therefore, $D_\l=
u_\l'K[h]$ for some element $u_\l' \in D_\l$. Then $u_\l = u_\l'p$
for some nonzero element $p\in K[h]$, and so $$u_\l'K[h]= D_\l =
K[h]u_\l = K[h]u_\l'p\subseteq D_\l p = u_\l ' K[h]p\subseteq
u_\l'K[h].$$
 This gives the equality $K[h]= K[h] p$, and so $p\in K^*$. This
 implies the equality $K[h]u_\l = u_\l H[h]$.

$\noindent $

{\em Claim 2: $x=\mu u_1$ and $y=\nu u_{-1}$ for some $\mu , \nu
\in K^*$.}

$\noindent $

By Claim 1, $x=u_1\mu$ and $y = \nu u_{-1}$ for some  non-zero
polynomials $\mu , \nu \in K[h]$. Then $h = yx= \nu u_{-1}
u_1\mu$. Since the polynomial $u_{-1}u_1\in K[h]$ is not  a scalar
polynomial (as the product of non-scalar elements in the Weyl
algebra $A_1$ is a non-scalar element), the polynomials $\nu $ and
$\mu$ must be  non-zero scalars since
$$ 1= \deg_h(h) = \deg_h( \l ) +\deg_h(\mu )
+\deg_h(u_{-1}u_1).$$ For every eigenvalue $\l \in \Ev (h)$, the
vector space $D_\l $ is invariant under the linear maps
$\d_y(\cdot ) x$ and $\d_x(\cdot ) y$.

\begin{lemma}\label{b23Dec10}
Let $d\in \{ \d_y(\cdot ) x, \d_x(\cdot ) y\}$, $v= v_{\rho ,
\eta}$ where $(\rho , \eta )\in \Z^2$ with $\rho + \eta >0$.  Then
$v(d(a^n))= v( a^{n-1}d(a))$ and $\D_{d, v}(a^n) = \D_{d,v}(a)$
for all elements $ a\in A_1 \backslash \ker (d)$ and $ n\geq 1$.
\end{lemma}

{\it Proof}. The second equality (i.e. $\D_{d, v}(a^n) =
\D_{d,v}(a)$) follows at once from the first one. In view of
existence of the automorphism of the Weyl algebra $A_1':= K\langle
x,y\rangle\ra A_1'$, $x\mapsto y$, $y\mapsto -x$, it suffices to
prove the first equality only  for $\d = \d_y(\cdot ) x$. Let us
prove by induction on $n$ that 
\begin{equation}\label{bana}
d(a^n) = na^{n-1}d(a)+\cdots \;\; {\rm for \; all}\;\; a\in
A_1\backslash \ker (\d_y)
\end{equation}
where the three dots denote smaller terms, i.e. $v(\cdots )
<v(na^{n-1}d(a))$. The claim is trivially true for $n=1$. So, let
$n>2$. By (\ref{dab}),
$$ d(a^n) = d(a\cdot a^{n-1}) = d(a) a^{n-1} +ad(a^{n-1}) + d(a) h^{-1}y[a^{n-1}, x]
.$$ Then, by (\ref{vrab}), 
\begin{equation}\label{bana1}
v(d(a) h^{-1}y[a^{n-1}, x])  \leq v( d(a)) - v(h) +v(y)
+v(a^{n-1})+v(x) - \rho -\eta = v(a^{n-1} d(a))-\rho -\eta .
\end{equation}
By induction on $n$, $a\cdot d(a^{n-1}) = (n-1) a\cdot a^{n-2}
d(a) +\cdots =(n-1)a^{n-1} d(a)+\cdots $. Then
\begin{eqnarray*}
d(a^n) &= & (a^{n-1} d(a)+\cdots )  +((n-1) a^{n-1} d(a) +\cdots )
+d(a) h^{-1}y[a^{n-1},x] \\
& =&(na^{n-1} d(a) +\cdots )+d(a)h^{-1}y[a^{n-1},x] = na^{n-1}
d(a) +\cdots . \;\; \Box
\end{eqnarray*}

\begin{lemma}\label{a1Jan11}
\begin{enumerate}
\item Each element of the basis $\{ y^n\cdot y^ix^i, y^ix^i,
y^ix^i\cdot x^n\, | \, n\geq 1, i\geq 0\}$ of the Weyl algebra
$A_1'$ is an eigenvector for the linear map $d= \d_y(\cdot ) x$:
$$ d(y^n\cdot y^ix^i)= iy^n\cdot y^ix^i, \;\; d(y^ix^i)= iy^ix^i, \;\;
d(y^ix^i\cdot x^n)=(i+n)y^ix^i\cdot x^n.$$ \item Each element of
the basis $\{ x^iy^i\cdot y^n, x^iy^i, x^n\cdot x^iy^i \, | \,
n\geq 1, i\geq 0\}$ of the Weyl algebra $A_1'$ is an eigenvector
for the linear map $d'= \d_x(\cdot ) y$:
$$ d'( x^iy^i\cdot y^n)= -(i+n)  x^iy^i\cdot y^n, \;\; d'(x^iy^i)= -ix^iy^i, \;\;
d'(x^n\cdot x^iy^i)=-i x^n\cdot x^iy^i.$$
\end{enumerate}
\end{lemma}

{\it Proof}.  Straightforward. $\;\;\Box $


\begin{corollary}\label{b1Jan11}
Let $d\in \{ \d_y(\cdot ) x, \d_x(\cdot ) y \}$, $v= v_{\rho ,
\eta}$ where $(\rho , \eta ) \in \Z^2$ with $\rho +\eta >0$, and
$\D = \D_{d , v}$. Then $v(d(a)) = v(a)$ for all $a\in
A_1'\backslash \ker (d)$, i.e. $\D (a) =0$ for all $a\in
A_1'\backslash \ker (d)$.
\end{corollary}

{\it Proof}. This follows at once from Lemma \ref{a1Jan11}. $\Box
$

$\noindent $

{\em Claim 3: $\th =1$, i.e. $n=1$}.

$\noindent $

Suppose that $n\neq 1$, we seek a contradiction. Fix a pair $(\rho
, \eta )\in \Z^2$ with $\rho + \eta >0$. Let $v = v_{\rho , \eta
}$, $d\in \{ \d_y(\cdot ) x, \d_x(\cdot ) y\}$, and $\D = \D_{d ,
v}$. Since $(u_\frac{1}{n})^n \in D_{n \frac{1}{n}} = D_1= xK[h]$
(by Claims 1 and 2), there exists  a nonzero polynomial $\alpha
\in K[H]\backslash \{ 0\}$ such that $(u_\frac{1}{n})^n= x\alpha$.
Clearly, $u_\frac{1}{n}\not\in \ker (d)$ since $\ker (\d_y(\cdot
)x) = \ker (\d_y) = K[y]$ and $\ker (\d_x(\cdot ) y) = \ker (\d_x)
= K[x]$, by Theorem \ref{GGV-Cen}, and $K[x]+K[y]\subseteq
\bigoplus_{i\in \Z } D_i$. By Lemma \ref{b23Dec10} and Corollary
\ref{b1Jan11}, $ \D (u_\frac{1}{n}) = \D ((u_\frac{1}{n})^n) = \D
( x\alpha ) =0$. Therefore, $d(u_\frac{1}{n}) = u_\frac{1}{n}\g_d$
for some $\g_d\in K^*$ since $d (D_\frac{1}{n})\subseteq
D_\frac{1}{n}$ and $D_\frac{1}{n}= u_\frac{1}{n}K[h]$ (Claim 1).
So, $ d(u_\frac{1}{n}) = u_\frac{1}{n}\g $ and $d' (u_\frac{1}{n})
= u_\frac{1}{n}\g'$ for some $\g , \g'\in K^*$ where $d =
\d_y(\cdot ) x$ and $d'= \d_x(\cdot ) y$. Let $\d = \d_x\d_y$.
Then $\d (D_\frac{1}{n}) \subseteq D_\frac{1}{n}$, and so $\d
(u_\frac{1}{n}) = u_\frac{1}{n}\alpha$ for some $\alpha \in K[h]$.
By (\ref{ddp}),
$$ u_\frac{1}{n}\g\g'= dd'(u_\frac{1}{n}) = \d
(u_\frac{1}{n})h=u_\frac{1}{n}\alpha h,$$ and so $\alpha \neq 0$
and $0=\deg_h( \g\g') = \deg_h(\alpha h) \geq 1$, a contradiction.

$\noindent $

{\em Claim 4: For all $n\geq 1$, $u_n= \l_nx^n$ and $u_{-n}=
\l_{-n}y^n$ for some scalars $\l_n, \l_{-n}\in K^*$.}

 $\noindent $

In view of existence of the automorphism of the Weyl algebra
$A_1'\ra A_1'$, $ x\mapsto y$, $y\mapsto -x$ $(h\mapsto -h+1)$, it
suffices to prove Claim 4 for positive $n$. By Claim 1, $x^n =
u_n\alpha$ for some element $ \alpha \in K[h]\backslash \{ 0\}$.
Suppose that $\alpha \not\in K^*$, we seek a contradiction. We
keep the notation of the proof of Claim 3. In addition, we assume
that $\rho >0$ and $\eta >0$. Then $d(u_n) = u_n \g $ and $d'(u_n)
= u_n \g'$ for some elements $\g , \g' \in K[h]\backslash \{ 0\}$,
by Claim 1 and since $\ker_{A_1}(d ) = \ker_{A_1}(\d_y) = K[y]$
and  $\ker_{A_1}(d' ) = \ker_{A_1}(\d_x) = K[x]$, by Theorem
\ref{GGV-Cen} (where $d:= \d_y (\cdot ) x$ and $ d' := \d_x(\cdot
) y$). Let $ m:= \eta n = v(x^n) = v( u_n \alpha )$. Applying the
linear map $ d$ to the equality $x^n = u_n\alpha$ and using
(\ref{dab}) we obtain the equality
$$nx^n = d(x^n) = d( u_n \alpha ) = d(u_n) \alpha +u_nd(\alpha ) -
d(u_n)h^{-1}yx (\alpha - \s^{-1} (\alpha )) = u_n \g \alpha + u_n
d(\alpha ) -u_n \g (\alpha - \s^{-1} (\alpha ))$$ where $\s :
h\mapsto h-1$ is an automorphism of the polynomial algebra $K[h]$.
Notice that $\deg_h( \alpha - \s^{-1} (\alpha )) <\deg_h(\alpha
)$, and so $v(u_n\g (\alpha - \s^{-1} (\alpha )))<v(u_n\g \alpha
)$. By the assumption, $\alpha \in K[h]\backslash K$, hence
$d(\alpha ) \neq 0$ since $\ker (d) = \ker (\d_y) = K[y]$ and
$K[h]\cap K[y]= K$. Then, by Corollary \ref{b1Jan11}, $v(d(\alpha
)) = v( \alpha )$, and so $v(u_nd(\alpha )) = v(u_n) +v(d(\alpha
)) = v(u_n) +v(\alpha )= v(u_n \alpha ) = m$. Therefore,
$$ m+(\rho +\eta ) \deg_h(\g ) = v( u_n \g \alpha )=
 v( u_n \g \alpha - u_n\g (\alpha - \s^{-1} (\alpha )) \leq \max \{
v(nx^n), v( u_nd(\alpha ))\} = m.$$ Therefore, $\deg_h(\g ) =0$
since $\rho +\eta >0$, and so $\g \in K^*$.

Similarly, applying the linear map $d' := \d_x(\cdot ) y$ to the
equality $x^n = u_n \alpha $ and using (\ref{sdab}) we obtain the
equality
\begin{eqnarray*}
 0&=& d'(u_n \alpha ) = d'(u_n) \alpha + u_n d'(\alpha ) +
 d'(u_n) (h-1)^{-1} x[\alpha , y]\\
 &=&
  u_n\g' \alpha + u_nd'(\alpha ) -u_n \g' (xy)^{-1} xy (\alpha - \s (\alpha ))  \\
 &=& u_n \g ' \alpha + u_n d'(\alpha ) -u_n \g ' (\alpha - \s
 (\alpha )).
\end{eqnarray*}
Notice that $\deg_h( \alpha - \s (\alpha ))<\deg_h(\alpha )$, and
so
$$v(\alpha - \s (\alpha )) = v(h) \deg_h(\alpha -
\s (\alpha ))<v(h) \deg_h(\alpha ) = v(\alpha ),$$ hence $v(u_n\g'
\alpha ) >v(u_n \g' (\alpha - \s (\alpha )))$. By the assumption,
$\alpha \in K[h]\backslash K$, hence $d'(\alpha ) \neq 0$ since
$\ker (d')= \ker (\d_x)=K[x]$ and $ K[h]\cap K[x]=K$. Then, by
Corollary \ref{b1Jan11}, $v(d'(\alpha )) = v(\alpha )$, and so
$$ v(u_n) +v(\g') +v(\alpha ) = v(u_n\g' \alpha ) = v(u_n
d'(\alpha ))=v(u_n)+v(d'(\alpha )) = v(u_n) +v(\alpha ).$$
Therefore, $0=v(\g') = v(h)\deg_h(\g ')$, i.e. $\g'\in K^*$.

 Recall that $\d = \d_x\d_y$.
Since $\d (D_n) \subseteq D_n$ and $D_n = u_n K[h]$, there exists
$\beta \in K[h]$ such that $\d (u_n) = u_n \beta$. By (\ref{ddp}),
$$ u_n \g\g' = dd'(u_n) = \d (u_n) h = u_n \beta h,$$
hence $\beta \neq 0$ and $ 0 = \deg_h(\g\g') = \deg_h(\beta h)
\geq 1$, a contradiction. The proof of Claim 4 is complete. By
Claims 1, 3 and 4, $D(h) = A_1'$. The proof of Theorem
\ref{22Dec10} is complete. $\Box$

\begin{proposition}\label{2Jan11}
Let $d= \d_y(\cdot ) x$, $d'= \d_x(\cdot ) y$, and $a\in A_1$. If
$(dd')^n (a) \in A_1'$ for some $n\geq 1$ then $a\in A_1'$.
\end{proposition}

{\it Proof}. We use induction on $n$ to prove the result. Notice
that $dd' = \d (\cdot ) h$ where $\d = \d_x\d_y$, see (\ref{ddp}).
If $n=1$, i.e. $\d (a) h \in A_1'$ then  $\d (a) \in \Frac (A_1')
\cap A_1= A_1'$ (Theorem \ref{X1Jan11}). Since $A_1'= N(\d , A_1)$
(Proposition \ref{A18Dec10}.(1)) and  $\d (a) \in A_1'$, we must
have $a\in N( \d , A_1) = A_1'$.

Suppose that $n>1$ and the statement is true for all $n'<n$. Then
$(dd')^n(a) = (dd')^{n-1} (dd' (a)) \in A_1'$, hence $dd'(a) \in
A_1'$ (by the inductive hypothesis), and so $a\in A_1'$, by the
case $n=1$.  $\Box $

$\noindent $

{\it Question 2. Let $a\in A_1'$. Is it true that $N(a, A_1') =
N(a, A_1)$?}

$\noindent $

{\it Question 3. Let $a\in A_1'$. Is it true that $D(a, A_1') =
D(a, A_1)$?}

Department of Pure Mathematics

University of Sheffield

Hicks Building

Sheffield S3 7RH

UK

email: v.bavula@sheffield.ac.uk

\end{document}